\documentclass[12pt]{article} 

\usepackage{latexsym, 
amscd, 
graphicx, amssymb}

\newtheorem{thm}{Theorem}[section]

\newcommand{\text}[1]{\mbox{\rm #1}}

\newcommand{\nn}{\nonumber \\}

\newcommand{\wt}{\mbox{\rm wt}\ }
\newcommand{\swt}{\mbox{\rm {\scriptsize wt}}\ }

\newcommand{\tr}{\mbox{\rm Tr}}
\newcommand{\A}{\mathcal{A}}
\newcommand{\Y}{\mathcal{Y}}
\newcommand{\C}{\mathbb{C}}
\newcommand{\Z}{\mathbb{Z}}

\newcommand{\Q}{\mathbb{Q}}
\newcommand{\N}{\mathbb{N}}

\title{ {\bf Vertex operator algebras, fusion rules and 
modular transformations} }
\date{}
\author{Yi-Zhi Huang}

\begin{document}

\bibliographystyle{alpha}
\maketitle

\begin{abstract} 
We discuss a recent proof by the author 
of a general version of the Verlinde conjecture 
in the framework of vertex operator algebras and 
the application of this result to the construction of 
modular tensor tensor category structure on the 
category of modules for a vertex operator algebra.
\end{abstract}

\renewcommand{\theequation}{\thesection.\arabic{equation}}
\renewcommand{\thethm}{\thesection.\arabic{thm}}
\setcounter{equation}{0}
\setcounter{thm}{0}
\setcounter{section}{-1}

\section{Introduction}

One of the most important discoveries by physicists in two-dimensional
conformal field theory is the famous relation between the fusion rules
and the action of the modular transformation $\tau\mapsto -1/\tau$ on
the space of vacuum characters. It states that this action of the
modular transformation diagonalizes the matrices formed by the fusion
rules.  This relation was first conjectured by Verlinde \cite{V} based
on a comparison between the fusion algebra of a rational conformal field
theory and an algebra arising from the study of the genus-one part of
the theory. Assuming the axioms for rational conformal field theories,
Moore and Seiberg \cite{MS1} proved this Verlinde conjecture by deriving
a fundamental set of polynomial equations.  Moore and Seiberg \cite{MS2}
also observed that the genus-zero part of these polynomial
equations is analogous to braided tensor categories.  The theory of
complete sets of such polynomial equations was then called ``modular
tensor category'' which was first suggested by I. Frenkel.
Later, this notion of modular tensor category was reformulated precisely
using the language of tensor categories (see, for example, 
\cite{T} and \cite{BK}
for expositions and references on modular tensor categories). 

Given a modular tensor category, we have fusion rules and an action of
modular transformations. It is not difficult to show 
(see, for example, \cite{BK}) that for the fusion rules
and the action of the modular transformation $\tau\mapsto -1/\tau$
given in this way, the
Verlinde conjecture holds. But this version of the Verlinde conjecture
is not the original one because the action of the modular transformation
is the one constructed from the modular tensor category, not the one on the
space of the vacuum characters of the corresponding conformal field
theory. The missing piece is an identification of 
the two actions of the same modular transformation $\tau\mapsto -1/\tau$,
or equivalently, is a mathematical construction of the modular
tensor category associated to a rational conformal field theory. 
Moreover, the starting point of 
the original work \cite{MS1} and \cite{MS2} 
of Moore and Seiberg is the axioms for rational conformal 
field theories, which are actually even more difficult to prove than 
the Verlinde conjecture. It is therefore desirable to formulate and prove
a general and mathematical version of the Verlinde conjecture which should 
give precise and natural conditions on the vertex operator algebras 
(chiral algebras) such that the Verlinde conjecture holds
for these algebras. 

Recently, the author was able to formulate and prove such a 
general and mathematical
version of the Verlinde conjecture. Using this result, the author
has also proved the rigidity and nondegeneracy of the 
semisimple braided tensor 
category structure constructed by Lepowsky and the author 
on the category of modules for the vertex operator 
algebra. In particular, 
the modular tensor category structure on the category of 
modules for the vertex operator algebra is mathematically constructed.
See \cite{H7} and \cite{H8} for details and see also 
\cite{H9} for an announcement of the results. 
In the present paper, we shall discuss this general and 
mathematical version 
of the Verlinde conjecture, its proof and its application to the 
proofs of the rigidity and nondegerneracy of the braided tensor 
category structure mentioned above.

\paragraph{Acknowledgment} I am grateful to J\"{u}rgen Fuchs for 
inviting me to give a talk at this conference.
This research is partially supported 
by NSF grant DMS-0401302.

\renewcommand{\theequation}{\thesection.\arabic{equation}}
\renewcommand{\thethm}{\thesection.\arabic{thm}}
\setcounter{equation}{0}
\setcounter{thm}{0}

\section{Vertex operator algebras and fusion rules}

Vertex (operator) algebras were introduced in 1986 by Borcherds 
\cite{B} in connection with representations of 
affine Lie algebras and the moonshine module for 
the Monster finite simple group, which was
conjectured to 
exist by Conway and Norton \cite{CN} and constructed using vertex operators
by Frenkel, Lepowsky and Meurman \cite{FLM1} \cite{FLM2}.
These algebras are essentially equivalent to 
chiral algebras (see, for example, \cite{MS2}) in physics, 
which were first studied systematically 
by Belavin, Polyakov and Zamolodchikov 
\cite{BPZ} in 1984, though without the name chiral algebra.
Here we explain briefly the basic concepts in the theory of 
vertex operator algebras. 

A vertex operator algebra is a graded vector space 
$V=\coprod_{n\in \mathbb{Z}}V_{(n)}$ equipped with
a vertex operator map $Y:V\otimes V\to V((z))$ (the space of 
formal Laurent series in $z$ with finitely many negative power terms), 
a vacuum $\mathbf{1}\in V$ and a conformal element $\omega\in V$.
These data must satisfy the following axioms which are very natural from
the point of view of conformal field theory in physics: One formulation of 
the main axioms is:
For $u_{1}, u_{2}, v\in V$, 
$v'\in V'=\coprod_{n\in \mathbb{Z}}V_{(n)}^{*}$, the series 
\begin{eqnarray*}
&\langle v', Y(u_{1}, z_{1})Y(u_{2}, z_{2})v\rangle&\\
&\langle v', Y(u_{2}, z_{2})Y(u_{1}, z_{1})v\rangle&\\
&\langle v', Y(Y(u_{1}, z_{1}-z_{2})u_{2}, z_{2})v\rangle&
\end{eqnarray*}
are absolutely convergent
in the regions $|z_{1}|>|z_{2}|>0$, $|z_{2}|>|z_{1}|>0$ and
$|z_{2}|>|z_{1}-z_{2}|>0$, respectively, to a common rational function 
in $z_{1}$ and $z_{2}$ with the only possible poles at $z_{1}, z_{2}=0$ and 
$z_{1}=z_{2}$. Other axioms include: 
$$\dim V_{(n)}<\infty$$
for $n\in \Z$,
$$V_{(n)}=0$$
when $n$ is sufficiently negative 
(these are called grading-restriction conditions);
$$Y(\mathbf{1}, z)=1,\;\;\;\;\lim_{z\to 0}Y(u, z)\mathbf{1}=u;$$
let $L(n): V\to V$ be defined by 
$Y(\omega, z)=\sum_{n\in \mathbb{Z}}L(n)z^{-n-2}$,
then 
$$[L(m), L(n)]=(m-n)L(m+n)+\frac{c}{12}(m^{3}-m)\delta_{m+n, 0}$$
($c$ is called the central charge of $V$),
$$\frac{d}{dz}Y(u, z)=Y(L(-1)u, z)$$
and
$$L(0)u=nu \;\;\;\;\;\;\mbox{\rm for}\;\;u\in V_{(n)}$$
($n$ is called the weight of $u$ and is denoted $\wt u$).
For $u\in V$, we write $Y(u, z)=\sum_{n\in \mathbb{Z}}u_{n}z^{-n-1}$.
Then for homogeneous $u\in V$, the maps $u_{n}: V\to V$ have weights
$\wt u-n-1$, that is, they map $V_{(m)}$ to $V_{(\swt u-n-1+m)}$.

For vertex operator algebras, we have the following important 
condition which was first introduced by Zhu \cite{Z}:

\vspace{1em}

\noindent {\bf $C_{2}$-cofiniteness condition:} Let $V$ be a vertex
operator algebra and $C_{2}(V)$ be the subspace of $V$ spanned by
elements of the form $u_{-2}v$ for $u, v\in V$.  Then we say that $V$ is
{\it $C_{2}$-cofinite} or satisfies the {\it $C_{2}$-cofiniteness
condition} if $\dim V/C_{2}(V)<\infty$. 

\vspace{1em}

In other words, the $C_{2}$-cofiniteness condition says that 
the coefficients of the terms in the first power of $z$ 
in $Y(u, z)v$ for all $u, v\in V$ span the whole space $V$ except for
a finite-dimensional subspace. 
This condition is an easy consequence of the 
condition that the spaces of genus-one 
conformal blocks are finite-dimensional. The author showed in 
\cite{H6} that together with some 
other conditions, the $C_{2}$-cofiniteness condition
also implies the finiteness of the dimensions of the spaces of genus-one 
conformal blocks. 
In practice, it is much easier to verify this condition than to
prove the finiteness of the dimensions of the spaces of 
genus-one conformal blocks.  For vertex operator algebras 
associated to the Wess-Zumino-Novikov-Witten
models, vertex operator algebras associated to minimal models, lattice
vertex operator algebras and the moonshine module vertex operator
algebra, the $C_{2}$-cofiniteness condition was stated to 
hold in Zhu's paper \cite{Z} and was verified by Dong-Li-Mason 
\cite{DLM} (see also \cite{AN} for the case of minimal models).

A $V$-module can be defined simply as a $\mathbb{C}$-graded vector space 
$W=\coprod_{n\in \mathbb{C}}W_{(n)}$ equipped with a vertex operator map 
$Y_{W}: V\otimes W\to W[[z, z^{-1}]]$ satisfying all those axioms for $V$
which still make sense. We also need the notion of 
$\mathbb{N}$-gradable weak $V$-module. An $\mathbb{N}$-gradable 
weak $V$-module 
is an $\mathbb{N}$-graded
vector space $W=\coprod_{n\in \mathbb{N}}W_{[n]}$ and a vertex operator map
\begin{eqnarray*}
Y: V\otimes W&\to &W[[z, z^{-1}]] \\
u\otimes w&\mapsto &Y(u, z)w=\sum_{n\in \mathbb{Z}}u_{n}z^{-n-1}
\end{eqnarray*}
satisfying all axioms for $V$-modules except that the 
condition $L(0)w=nw$ for $w\in 
W_{(n)}$ is replaced by $u_{k}w\in W_{[m-k-1+n]}$ for 
$u\in V_{(m)}$ and $w\in W_{[n]}$. 

For $V$-modules $W_{1}$, $W_{2}$ and $W_{3}$, an
intertwining operator of type ${W_{3}\choose W_{1}W_{2}}$
is a linear map $\mathcal{Y}: W_{1}\otimes
W_{2}\to W_{3}\{z\}$, where 
$W_{3}\{z\}$ is the space of all series in complex powers of $z$ with
coefficients in $W_{3}$,  satisfying all those axioms for $V$ which still 
make sense. For example, for $u\in V$, $w_{1}\in W_{1}$, $w_{2}\in W_{2}$
and $w'_{3}\in W'_{3}=\coprod_{n\in \mathbb{C}}(W_{3})_{(n)}^{*}$,
\begin{eqnarray*}
&\langle w'_{3}, Y_{W_{3}}(u, z_{1})\mathcal{Y}(w_{1}, z_{2})w_{2}\rangle&\\
&\langle w'_{3}, \mathcal{Y}(w_{1}, z_{2})Y_{W_{2}}(u, z_{1})w_{2}\rangle&\\
&\langle w'_{3}, \mathcal{Y}(Y_{W_{1}}(u, z_{1}-z_{2})w_{1}, z_{2})
w_{2}\rangle&
\end{eqnarray*}
are absolutely convergent
in the regions $|z_{1}|>|z_{2}|>0$, $|z_{2}|>|z_{1}|>0$ and
$|z_{2}|>|z_{1}-z_{2}|>0$, respectively, to a common (multivalued)
analytic function 
in $z_{1}$ and $z_{2}$ with the only possible singularities
(branch points)
at $z_{1}, z_{2}=0$ and 
$z_{1}=z_{2}$. Also
$$\frac{d}{dz}\mathcal{Y}(w_{1}, z)=\Y(L(-1)w_{1}, z).$$
We denote the space of intertwining operator of type 
${W_{3}\choose W_{1}W_{2}}$ by $\mathcal{V}_{W_{1}W_{2}}^{W_{3}}$.
The dimension of $\mathcal{V}_{W_{1}W_{2}}^{W_{3}}$ is the 
fusion rule $N_{W_{1}W_{2}}^{W_{3}}$. 

\renewcommand{\theequation}{\thesection.\arabic{equation}}
\renewcommand{\thethm}{\thesection.\arabic{thm}}
\setcounter{equation}{0}
\setcounter{thm}{0}

\section{Fusing and braiding isomorphisms}

In the definitions of vertex operator algebra,
module and intertwining operator, we see that the 
main axioms are all about products and iterates of 
vertex operators or products and iterates of 
one vertex operator and one intertwining operator.
It is natural to expect that products and iterates of intertwining 
operators should have similar properties. Indeed, 
it was proved by the author in \cite{H1} and \cite{H5}
that for vertex operator algebras satisfying suitable 
finiteness and reductivity conditions, intertwining operators
satisfy associativity and commutativity properties. 
These properties give fusing and braiding isomorphisms.

Note that the associativity for intertwining operators is a strong
version of the operator product expansion of
``chiral vertex operators'' (which is equivalent to intertwining
operators for vertex operator algebras). 
In fact, in the important work
\cite{MS1} and \cite{MS2} of Moore and Seiberg, in addition to other axioms for 
rational conformal field theories, the operator product expansion of
chiral vertex operators was a basic assumption.
In physics, though there were calculations for particular
operators in particular examples, operator product expansion was
not a mathematical theorem and calculations based on it are by no
means simple. The
associativity for intertwining operators proved in \cite{H1} and \cite{H5}
under suitable conditions is in
fact a strong version of operator product expansion for
intertwining operators in the following sense: 
The associativity states that the
product of two intertwining operators is equal to the iterate of two
other intertwining operators in a suitable region. If we expand the
intertwining operator inside the other intertwining operator in the
iterate, we immediately obtain the operator product expansion of
intertwining operators.  On the other hand, 
in \cite{H1}, the author proved associativity
for intertwining operators under the assumption that the convergence
and extension property and some other conditions hold.  Operator
product expansion for chiral vertex operators (intertwining operators)
assumed in \cite{MS1} and \cite{MS2}
actually implies the convergence and extension
property. So the operator product expansion together with some other
conditions implies associativity for intertwining operators. Since the
associativity immediately gives the operator product expansion, 
but on the other hand,
only together 
with some other conditions the operator product expansion gives
the associativity, we see that the associativity is indeed stronger
than the operator product expansion.

Let $V$ be a simple vertex operator algebra of central charge $c$ and 
$V'$ the contragredient 
module of $V$ as a $V$-module. In the remaining part
of the present paper, we
shall always assume that $V$ satisfies the following conditions:

\begin{enumerate}

\item $V_{(n)}=0$ for $n<0$, $V_{(0)}=\mathbb{C}\mathbf{1}$ and 
$V'$ is isomorphic to $V$.

\item Every $\mathbb{N}$-gradable weak $V$-module is completely 
reducible.

\item $V$ is $C_{2}$-cofinite, that is, $\dim V/C_{2}(V)<\infty$. 

\end{enumerate}

Note that finitely generated $\N$-gradable weak $V$-modules are 
what naturally appear in the proofs of 
the theorems on genus-zero and genus-one 
correlation functions. Thus Condition 2 is natural and necessary because 
the Verlinde conjecture concerns $V$-modules,
not finitely generated $\N$-gradable weak $V$-modules. 
For vertex operator algebras 
associated to affine Lie algebras (Wess-Zumino-Novikov-Witten
models) and vertex operator algebras associated to the Virasoro algebra
(minimal models), Condition 2 can be verified easily by reformulating 
the corresponding complete reducibility results in terms of 
the representation theory
of affine Lie algebras and the Virasoro algebra. 

From \cite{DLM},  we know that there are only finitely many
inequivalent irreducible $V$-modules.
Let $\mathcal{A}$ be the set of
equivalence classes of irreducible $V$-modules. We denote the
equivalence class containing $V$ by $e$. For each $a\in \mathcal{A}$,
we choose a representative $W^{a}$ of $a$. Note that the contragredient
module of an irreducible module is also irreducible (see \cite{FHL}).
So we have a map
\begin{eqnarray*}
^{\prime}:
\mathcal{A}&\to& \mathcal{A}\\
a&\mapsto& a'.
\end{eqnarray*}
From \cite{AM} and \cite{DLM}, we know that irreducible $V$-modules 
are in fact graded by rational numbers. 
Thus for $a\in \A$, 
there exist $h_{a}\in \Q$ such that 
$W^{a}=\coprod_{n\in h_{a}+\N}W^{a}_{(n)}$.

Let $\mathcal{V}_{a_{1}a_{2}}^{a_{3}}$ for $a_{1}, a_{2}, a_{3}
\in \mathcal{A}$
be the space of intertwining operators of 
type ${W^{a_{3}}\choose W^{a_{1}}W^{a_{2}}}$ and 
$N_{a_{1}a_{2}}^{a_{3}}$ for $a_{1}, a_{2}, a_{3}
\in \mathcal{A}$ the fusion rule, that is, the dimension of 
the space of intertwining operators of 
type ${W^{a_{3}}\choose W^{a_{1}}W^{a_{2}}}$. For any 
$\Y\in \mathcal{V}_{a_{1}a_{2}}^{a_{3}}$, we know from \cite{FHL}
that 
for $w_{a_{1}}\in W^{a_{1}}$ and $w_{a_{2}}\in W^{a_{2}}$
\begin{equation}\label{y-delta}
\Y(w_{a_{1}}, x)w_{a_{2}}\in x^{\Delta(\Y)}W^{a_{3}}[[x, x^{-1}]],
\end{equation}
where 
$$\Delta(\mathcal{Y})=h_{a_{3}}-h_{a_{1}}-h_{a_{2}}.$$

From \cite{GN},
\cite{L},
\cite{AN},
\cite{H6}, we also know that the fusion rules 
$N_{a_{1}a_{2}}^{a_{3}}$ for 
$a_{1}, a_{2}, a_{3}\in \mathcal{A}$ are all finite. 
For $a\in \A$, let $\mathcal{N}(a)$ be the matrix 
whose entries are $N_{aa_{1}}^{a_{2}}$ for $a_{1}, a_{2}\in \A$,
that is, 
$$\mathcal{N}(a)=(N_{aa_{1}}^{a_{2}}).$$

We also need matrix elements of fusing and braiding isomorphisms. 
In the proof of the Verlinde conjecture, we need 
to use several bases of one space of intertwining 
operators. We shall  use $p=1, 2, 3, 4, 5, 6, \dots$
to label different bases. 
For $p=1, 2, 3, 4, 5, 6, \dots$ and $a_{1}, a_{2}, 
a_{3}\in \mathcal{A}$, let 
$\{\Y_{a_{1}a_{2}; i}^{a_{3}; (p)}\;|\; i=1, \dots, 
N_{a_{1}a_{2}}^{a_{3}}\}$, be a
basis of $\mathcal{V}_{a_{1}a_{2}}^{a_{3}}$. 
For $a_{1}, \dots, a_{6}\in \mathcal{A}$, 
$w_{a_{1}}\in W^{a_{1}}$, $w_{a_{2}}\in W^{a_{2}}$,
$w_{a_{3}}\in W^{a_{3}}$, and $w'_{a_{4}}\in (W^{a_{4}})'$,
using the differential equations satisfied by the series
$$\langle w'_{a_{4}}, \Y_{a_{1}a_{5}; i}^{a_{4}; (1)}(w_{a_{1}}, x_{1})
\Y_{a_{2}a_{3}; j}^{a_{5}; (2)}(w_{a_{2}}, x_{2})w_{a_{3}}\rangle
|_{x_{1}^{n}=e^{n\log z_{1}}, \; x_{2}^{n}=e^{n\log z_{2}}, \;
n\in \Q}$$
and 
$$\langle w'_{a_{4}}, 
\Y_{a_{6}a_{3}; k}^{a_{4}; (3)}
(\Y_{a_{1}a_{2}; l}^{a_{6}; (4)}(w_{a_{1}}, x_{0})
w_{a_{2}}, x_{2})w_{a_{3}}\rangle
|_{x_{0}^{n}=e^{n\log (z_{1}-z_{2})}, \; x_{2}^{n}=e^{n\log z_{2}}, \;
n\in \Q},$$
it was proved in \cite{H5} that these series are convergent 
in the regions $|z_{1}|>|z_{2}|>0$ and $|z_{2}|>|z_{1}-z_{2}|>0$,
respectively. Note that for any $a_{1}$, $a_{2}$, $a_{3}$, $a_{4}$, $a_{5}$,
$a_{6}\in \A$, 
$\Y_{a_{1}a_{5}; i}^{a_{4}; (1)}\otimes \Y_{a_{2}a_{3}; j}^{a_{5}; (2)}$
and 
$\Y_{a_{6}a_{3}; l}^{a_{4}; (3)}
\otimes \Y_{a_{1}a_{2}; k}^{a_{6}; (4)}$
are bases of $\mathcal{V}_{a_{1}a_{5}}^{a_{4}}\otimes 
\mathcal{V}_{a_{2}a_{3}}^{a_{5}}$ and 
$\mathcal{V}_{a_{6}a_{3}}^{a_{4}}
\otimes \mathcal{V}_{a_{1}a_{2}}^{a_{6}}$, respectively. 
The associativity of intertwining operators proved and 
studied in \cite{H1}, \cite{H4} and
\cite{H5} says that 
there exist 
$$F(\Y_{a_{1}a_{5}; i}^{a_{4}; (1)}\otimes \Y_{a_{2}a_{3}; j}^{a_{5}; (2)}; 
\Y_{a_{6}a_{3}; l}^{a_{4}; (3)}
\otimes \Y_{a_{1}a_{2}; k}^{a_{6}; (4)}) \in \C$$
for $a_{1}, \dots, a_{6}\in \mathcal{A}$, $i=1, \dots, 
N_{a_{1}a_{5}}^{a_{4}}$, $j=1, \dots, 
N_{a_{2}a_{3}}^{a_{5}}$, $k=1, \dots, 
N_{a_{6}a_{3}}^{a_{4}}$, $l=1, \dots, 
N_{a_{1}a_{2}}^{a_{6}}$
such that 
\begin{eqnarray}\label{assoc}
\lefteqn{\langle w'_{a_{4}}, \Y_{a_{1}a_{5}; i}^{a_{4}; (1)}(w_{a_{1}}, x_{1})
\Y_{a_{2}a_{3}; j}^{a_{5}; (2)}(w_{a_{2}}, z_{2})w_{a_{3}}\rangle
|_{x_{1}^{n}=e^{n\log z_{1}}, \; x_{2}^{n}=e^{n\log z_{2}}, \;
n\in \Q}}\nn
&&=\sum_{a_{6}\in \A}
\sum_{k=1}^{N_{a_{6}a_{3}}^{a_{4}}}\sum_{l=1}^{N_{a_{1}a_{2}}^{a_{6}}}
F(\Y_{a_{1}a_{5}; i}^{a_{4}; (1)}\otimes \Y_{a_{2}a_{3}; j}^{a_{5}; (2)}; 
\Y_{a_{6}a_{3}; l}^{a_{4}; (3)}
\otimes \Y_{a_{1}a_{2}; k}^{a_{6}; (4)})\cdot\nn
&&\quad \cdot 
\langle w'_{a_{4}}, 
\Y_{a_{6}a_{3}; k}^{a_{4}; (3)}(\Y_{a_{1}a_{2}; l}^{a_{6}; (4)}(w_{a_{1}}, z_{1}-z_{2})
w_{a_{2}}, z_{2})w_{a_{3}}\rangle
|_{x_{0}^{n}=e^{n\log (z_{1}-z_{2})}, \; x_{2}^{n}=e^{n\log z_{2}}, \;
n\in \Q}\nn
&&
\end{eqnarray}
when $|z_{1}|>|z_{2}|>|z_{1}-z_{2}|>0$, 
for $a_{1}, \dots, a_{5}\in \A$, $w_{a_{1}}\in W^{a_{1}}$, 
$w_{a_{2}}\in W^{a_{2}}$, $w_{a_{3}}\in W^{a_{3}}$,  
$w'_{a_{4}}\in (W^{a_{4}})'$, $i=1, \dots, 
N_{a_{1}a_{5}}^{a_{4}}$ and $j=1, \dots, 
N_{a_{2}a_{3}}^{a_{5}}$. The numbers 
$$F(\Y_{a_{1}a_{5}; i}^{a_{4}; (1)}\otimes \Y_{a_{2}a_{3}; j}^{a_{5}; (2)}; 
\Y_{a_{6}a_{3}; k}^{a_{4}; (3)}
\otimes \Y_{a_{1}a_{2}; l}^{a_{6}; (4)})$$
are the matrix elements of the fusing isomorphism, that is, 
these numbers together give a matrix which represents a linear 
isomorphism, called fusing isomorphism, from 
$$\coprod_{a_{1}, a_{2}, a_{3}, a_{4}, a_{5}\in \A}
\mathcal{V}_{a_{1}a_{5}}^{a_{4}}\otimes 
\mathcal{V}_{a_{2}a_{3}}^{a_{5}}$$
to
$$\coprod_{a_{1}, a_{2}, a_{3}, a_{4}, a_{6}\in \A}
\mathcal{V}_{a_{6}a_{3}}^{a_{4}}
\otimes \mathcal{V}_{a_{1}a_{2}}^{a_{6}}.$$

By the commutativity of intertwining operators proved 
and studied in \cite{H2},
\cite{H4}
and \cite{H5}, 
there exist 
$$B^{(r)}(\Y_{a_{1}a_{5}; i}^{a_{4}; (1)}
\otimes \Y_{a_{2}a_{3}; j}^{a_{5}; (2)}; 
\Y_{a_{2}a_{6}; l}^{a_{4}; (3)}
\otimes \Y_{a_{1}a_{3}; k}^{a_{6}; (4)}) \in \C$$
for $r\in \Z$, $a_{1}, \dots, a_{6}\in \mathcal{A}$, $i=1, \dots, 
N_{a_{1}a_{5}}^{a_{4}}$, $j=1, \dots, 
N_{a_{2}a_{3}}^{a_{5}}$, $k=1, \dots, 
N_{a_{2}a_{6}}^{a_{4}}$, $l=1, \dots, 
N_{a_{1}a_{3}}^{a_{6}}$, 
such that 
the analytic extension of the single-valued analytic function
$$\langle w'_{a_{4}}, \Y_{a_{1}a_{5}; i}^{a_{4}; (1)}(w_{a_{1}}, x_{1})
\Y_{a_{2}a_{3}; j}^{a_{5}; (2)}(w_{a_{2}}, x_{2})w_{a_{3}}\rangle
|_{x_{1}^{n}=e^{n\log z_{1}}, \; x_{2}^{n}=e^{n\log z_{2}}, \;
n\in \Q}$$
on the region $|z_{1}|>|z_{2}|>0$, $0\le \arg z_{1}, \arg z_{2}<2\pi$
along the path 
$$t \mapsto \left(\frac{3}{2}
-\frac{e^{(2r+1)\pi i t}}{2}, \frac{3}{2}
+\frac{e^{(2r+1)\pi i t}}{2}\right)$$ 
to the region $|z_{2}|>|z_{1}|>0$, $0\le \arg z_{1}, \arg z_{2}<2\pi$
is 
\begin{eqnarray*}
\lefteqn{\sum_{a_{6}\in \A}
\sum_{k=1}^{N_{a_{2}a_{6}}^{a_{4}}}\sum_{l=1}^{N_{a_{1}a_{3}}^{a_{6}}}
B^{(r)}(\Y_{a_{1}a_{5}; i}^{a_{4}; (1)}\otimes \Y_{a_{2}a_{3}; j}^{a_{5}; (2)};
\Y_{a_{2}a_{6}; k}^{a_{4}; (3)}\otimes 
\Y_{a_{1}a_{3}; l}^{a_{6}; (4)})\cdot}\nn
&&\quad\quad\quad\quad\cdot 
\langle w'_{a_{4}}, 
\Y_{a_{2}a_{6}; k}^{a_{4}; (3)}(
w_{a_{2}}, z_{1})\Y_{a_{1}a_{3}; l}^{a_{6}; (4)}(w_{a_{1}}, z_{2})
w_{a_{3}}\rangle|_{x_{1}^{n}=e^{n\log z_{1}}, \; x_{2}^{n}=e^{n\log z_{2}}, \;
n\in \Q}.
\end{eqnarray*}
The numbers 
$$B^{(r)}(\Y_{a_{1}a_{5}; i}^{a_{4}; (1)}\otimes \Y_{a_{2}a_{3}; j}^{a_{5}; (2)};
\Y_{a_{2}a_{6}; k}^{a_{4}; (3)}\otimes \Y_{a_{1}a_{3}; l}^{a_{6}; (4)})$$
are the matrix elements of the braiding isomorphism.

We need an action of $S_{3}$ on the space 
$$\mathcal{V}=\coprod_{a_{1}, a_{2}, a_{3}\in 
\mathcal{A}}\mathcal{V}_{a_{1}a_{2}}^{a_{3}}.$$
For $a_{1}, a_{2}, a_{3}\in \mathcal{A}$, we have isomorphisms
$\Omega_{-r}: \mathcal{V}_{a_{1}a_{2}}^{a_{3}} \to 
\mathcal{V}_{a_{2}a_{1}}^{a_{3}}$ and 
$A_{-r}: \mathcal{V}_{a_{1}a_{2}}^{a_{3}} \to 
\mathcal{V}_{a_{1}a'_{3}}^{a'_{2}}$ for $r\in \Z$ (see \cite{HL2}).
For $a_{1}, a_{2}, a_{3}\in \A$,
$\mathcal{Y}\in \mathcal{V}_{a_{1}a_{2}}^{a_{3}}$, 
we define 
\begin{eqnarray*}
\sigma_{12}(\mathcal{Y})&=&e^{\pi i \Delta(\mathcal{Y})}
\Omega_{-1}(\mathcal{Y})\\
&=&e^{-\pi i \Delta(\mathcal{Y})}\Omega_{0}(\mathcal{Y}),\\
\sigma_{23}(\mathcal{Y})&=&e^{\pi i h_{a_{1}}}
A_{-1}(\mathcal{Y})\\
&=&e^{-\pi i h_{a_{1}}}A_{0}(\mathcal{Y}).
\end{eqnarray*}
Then the actions $\sigma_{12}$ and 
$\sigma_{23}$ of $(12)$ and 
$(23)$ on $\mathcal{V}$ defined above
generate a left action of $S_{3}$ on $\mathcal{V}$.

We now  choose
a basis $\mathcal{Y}_{a_{1}a_{2}; i}^{a_{3}}$, $i=1, \dots, 
N_{a_{1}a_{2}}^{a_{3}}$, 
of $\mathcal{V}_{a_{1}a_{2}}^{a_{3}}$ for each triple 
$a_{1}, a_{2}, a_{3}\in \mathcal{A}$. 
For $a\in \A$, we choose $\Y_{ea; 1}^{a}$ to be the vertex operator 
$Y_{W^{a}}$ defining the module structure on $W^{a}$ and we choose 
$\Y_{ae; 1}^{a}$ to be the intertwining operator defined using 
the action of $\sigma_{12}$,
\begin{eqnarray*}
\Y_{ae; 1}^{a}(w_{a}, x)u&=&\sigma_{12}(\Y_{ea; 1}^{a})(w_{a}, x)u\nn
&=&e^{xL(-1)}\Y_{ea; 1}^{a}(u, -x)w_{a}\nn
&=&e^{xL(-1)}Y_{W^{a}}(u, -x)w_{a}
\end{eqnarray*}
for $u\in V$ and $w_{a}\in W^{a}$. 
Since $V'$ as a $V$-module is isomorphic to $V$, we have 
$e'=e$. From \cite{FHL}, we know that there is a nondegenerate
invariant 
bilinear form $(\cdot, \cdot)$ on $V$ such that $(\mathbf{1}, 
\mathbf{1})=1$. 
We choose $\Y_{aa'; 1}^{e}=\Y_{aa'; 1}^{e'}$
to be the intertwining operator defined using the action of 
$\sigma_{23}$ by
$$\Y_{aa'; 1}^{e'}=\sigma_{23}(\Y_{ae; 1}^{a}),$$
that is,
$$(u, \Y_{aa'; 1}^{e'}(w_{a}, x)w'_{a})
=e^{\pi i h_{a}}\langle \Y_{ae; 1}^{a}(e^{xL(1)}(e^{-\pi i}x^{-2})^{L(0)}w_{a}, x^{-1})u, 
w'_{a}\rangle$$
for $u\in V$, $w_{a}\in W^{a}$ and $w'_{a'}\in (W^{a})'$. Since the actions of
$\sigma_{12}$
and $\sigma_{23}$ generate the action of $S_{3}$ on $\mathcal{V}$, we have
$$\Y_{a'a; 1}^{e}=\sigma_{12}(\Y_{aa'; 1}^{e})$$
for any $a\in \mathcal{A}$.
When $a_{1}, a_{2}, a_{3}\ne e$, we choose 
$\mathcal{Y}_{a_{1}a_{2}; i}^{a_{3}}$, $i=1, \dots, 
N_{a_{1}a_{2}}^{a_{3}}$, to be an arbitrary basis
of $\mathcal{V}_{a_{1}a_{2}}^{a_{3}}$.
Note that for each element 
$\sigma\in S_{3}$, 
$\{\sigma(\mathcal{Y})_{a_{1}a_{2}; i}^{a_{3}}\;|\;i=1, \dots, 
N_{a_{1}a_{2}}^{a_{3}}\}$, is also 
a basis of $\mathcal{V}_{a_{1}a_{2}}^{a_{3}}$.

\renewcommand{\theequation}{\thesection.\arabic{equation}}
\renewcommand{\thethm}{\thesection.\arabic{thm}}
\setcounter{equation}{0}
\setcounter{thm}{0}

\section{Modular Invariance}

We discuss modular invariance briefly in this section.
Let $q_{\tau}=e^{2\pi i\tau}$
for $\tau \in \mathbb{H}$ ($\mathbb{H}$ is the upper-half plane). 
We consider the $q_{\tau}$-traces of the vertex operators 
$Y_{W^{a}}$ for $a\in \mathcal{A}$ on 
the irreducible $V$-modules $W^{a}$ of the following form:
\begin{equation}\label{1-trace}
\tr_{W^{a}}
Y_{W^{a}}(e^{2\pi iz L(0)}u, e^{2\pi i z})
q_{\tau}^{L(0)-\frac{c}{24}}
\end{equation}
for $u\in V$ (recall that $c$ is the central charge of $V$). 
In \cite{Z}, under some conditions slightly different 
from (mostly stronger than) those we assume in this paper, 
Zhu proved that these $q$-traces are independent of 
$z$, are absolutely convergent
when $0<|q_{\tau}|<1$ and can be analytically extended to 
analytic functions of $\tau$ in the upper-half plane. 
We shall denote the analytic extension of (\ref{1-trace})
by 
$$E(\tr_{W^{a}}
Y_{W^{a}}(e^{2\pi iz L(0)}u, e^{2\pi i z})
q_{\tau}^{L(0)-\frac{c}{24}}).$$
In \cite{Z}, under his conditions mentioned above,
Zhu also proved the following modular invariance property:
For 
$$\left(\begin{array}{cc}a&b\\ c&d\end{array}\right)\in SL(2, \Z),$$
let $\tau'=\frac{a\tau+b}{c\tau+d}$. Then there exist 
unique $A_{a_{1}}^{a_{2}}\in \C$ for $a_{1}, a_{2}\in \mathcal{A}$
such that 
\begin{eqnarray*}
\lefteqn{E\left(\tr_{W^{a_{1}}}
Y_{W^{a_{1}}}\left(e^{\frac{2\pi iz}{c\tau+d}L(0)}
\left(\frac{1}{c\tau+d}\right)^{L(0)}
u, e^{\frac{2\pi i z}{c\tau+d}}\right)
q_{\tau'}^{L(0)-\frac{c}{24}}\right)}\nn
&&=\sum_{a_{2}\in \A}A_{a_{1}}^{a_{2}}
E(\tr_{W^{a_{2}}}
Y_{W^{a_{2}}}(e^{2\pi iz L(0)}u, e^{2\pi i z})
q_{\tau}^{L(0)-\frac{c}{24}})
\end{eqnarray*}
for $u\in V$. In \cite{DLM}, Dong, Li and Mason, 
among many other things,
improved Zhu's results above by showing that the results 
of Zhu above also hold for vertex operator algebras satisfying 
the conditions (slightly weaker than what) we assume in this paper. 
In particular, for 
$$\left(\begin{array}{cc}0&1\\ -1&0\end{array}\right)\in SL(2, \Z),$$
there exist unique $S_{a_{1}}^{a_{2}}\in \C$ for $a_{1}\in \mathcal{A}$
such that 
\begin{eqnarray*}
\lefteqn{E\left(\tr_{W^{a_{1}}}
Y_{W^{a_{1}}}\left(e^{-\frac{2\pi iz}{\tau}L(0)}
\left(-\frac{1}{\tau}\right)^{L(0)}
u, e^{-\frac{2\pi i z}{\tau}}\right)
q_{-\frac{1}{\tau}}^{L(0)-\frac{c}{24}}\right)}\nn
&&=\sum_{a_{2}\in \A}S_{a_{1}}^{a_{2}}
E(\tr_{W^{a_{2}}}
Y_{W^{a_{2}}}(e^{2\pi iz L(0)}u, e^{2\pi i z})
q_{\tau}^{L(0)-\frac{c}{24}})
\end{eqnarray*}
for $u\in V$. When $u=\mathbf{1}$, we see that the matrix
$S=(S_{a_{1}}^{a_{2}})$ actually acts on the space 
of spanned by the 
vacuum characters $\tr_{W^{a}}q_{\tau}^{L(0)-\frac{c}{24}}$
for $a\in \mathcal{A}$.

\renewcommand{\theequation}{\thesection.\arabic{equation}}
\renewcommand{\thethm}{\thesection.\arabic{thm}}
\setcounter{equation}{0}
\setcounter{thm}{0}

\section{The Verlinde conjecture and consequences}

In \cite{H7}, the author proved the following general version of 
the Verlinde conjecture in the framework of vertex operator algebras
(cf. Section 3 in \cite{V} and Section 4 in 
\cite{MS1}):

\begin{thm}\label{main}
Let $V$ be a vertex operator algebra satisfying the 
following conditions:

\begin{enumerate}

\item $V_{(n)}=0$ for $n<0$, $V_{(0)}=\mathbb{C}\mathbf{1}$ and 
$V'$ is isomorphic to $V$ as a $V$-module.

\item Every $\mathbb{N}$-gradable weak $V$-module is completely 
reducible.

\item $V$ is $C_{2}$-cofinite, that is, $\dim V/C_{2}(V)<\infty$. 

\end{enumerate}
Then for $a\in \A$, 
$$F(\Y_{ae; 1}^{a} \otimes \Y_{a'a; 1}^{e};
\Y_{ea; 1}^{a}\otimes \Y_{aa'; 1}^{e})\ne 0$$
and 
\begin{equation}\label{diag2}
\sum_{a_{1}, a_{3}\in \mathcal{A}}(S^{-1})_{a_{4}}^{a_{1}}
N_{a_{1}a_{2}}^{a_{3}}
S_{a_{3}}^{a_{5}}=\delta_{a_{4}}^{a_{5}}
\frac{(B^{(-1)})^{2}(\Y_{a_{4}e; 1}^{a_{4}}\otimes \Y_{a'_{2}a_{2}; 1}^{e};
\Y_{a_{4}e; 1}^{a_{4}}\otimes \Y_{a'_{2}a_{2}; 1}^{e})}
{F(\Y_{a_{2}e; 1}^{a_{2}} \otimes \Y_{a'_{2}a_{2}; 1}^{e};
\Y_{ea_{2}; 1}^{a_{2}}\otimes \Y_{a_{2}a'_{2}; 1}^{e})}.
\end{equation}
In particular, the matrix $S$ diagonalizes the
matrices $\mathcal{N}(a_{2})$ for all  $a_{2}\in \A$.
\end{thm}

The main work in the proof of this theorem is 
to prove the following formulas derived by
Moore and Seiberg \cite{MS1} from 
the axioms of rational conformal field
theories: For $a_{1}, a_{2}, a_{3}\in \A$, 
\begin{eqnarray}\label{formula1-2}
\lefteqn{\sum_{i=1}^{N_{a_{1}a_{2}}^{a_{3}}}\sum_{k=1}^{N_{a'_{1}a_{3}}^{a_{2}}}
F(\Y_{a_{2}e; 1}^{a_{2}}\otimes \Y_{a'_{3}a_{3}; 1}^{e}; 
\Y_{a'_{1}a_{3}; k}^{a_{2}}\otimes \Y_{a_{2}a'_{3}; i}^{a'_{1}})\cdot}\nn
&&\quad\quad\quad\quad\quad\cdot 
F(\Y_{a'_{1}a_{3}; k}^{a_{2}}\otimes \sigma_{123}(\Y_{a_{2}a'_{3}; i}^{a'_{1}});
\Y_{ea_{2}; 1}^{a_{2}}\otimes \Y_{a'_{1}a_{1}; 1}^{e})\nn
&&=N_{a_{1}a_{2}}^{a_{3}}
F(\Y_{a_{2}e; 1}^{a_{2}} \otimes \Y_{a'_{2}a_{2}; 1}^{e};
\Y_{ea_{2}; 1}^{a_{2}}\otimes \Y_{a_{2}a'_{2}; 1}^{e})
\end{eqnarray}
and
\begin{eqnarray}\label{formula2}
\lefteqn{\sum_{a_{4}\in \mathcal{A}}S_{a_{1}}^{a_{4}}
(B^{(-1)})^{2}(\Y_{a_{4}e; 1}^{a_{4}}\otimes \Y_{a'_{2}a_{2}; 1}^{e};
\Y_{a_{4}e; 1}^{a_{4}}\otimes \Y_{a'_{2}a_{2}; 1}^{e})
(S^{-1})_{a_{4}}^{a_{3}}}\nn
&&=\sum_{i=1}^{N_{a_{1}a_{2}}^{a_{3}}}\sum_{k=1}^{N_{a'_{1}a_{3}}^{a_{2}}}
F(\Y_{a_{2}e; 1}^{a_{2}}\otimes \Y_{a'_{3}a_{3}; 1}^{e}; 
\Y_{a'_{1}a_{3}; k}^{a_{2}}\otimes \Y_{a_{2}a'_{3}; i}^{a'_{1}})\cdot\nn
&&\quad\quad\quad\quad\quad\cdot 
F(\Y_{a'_{1}a_{3}; k}^{a_{2}}\otimes \sigma_{123}(\Y_{a_{2}a'_{3}; i}^{a'_{1}});
\Y_{ea_{2}; 1}^{a_{2}}\otimes \Y_{a'_{1}a_{1}; 1}^{e}).
\end{eqnarray}

The proof of the first formula (\ref{formula1-2})
uses mainly the works of Lepowsky and the author \cite{HL1}--\cite{HL4}
and of 
the author \cite{H1} \cite{H2} \cite{H4} and 
\cite{H5} on the tensor product theory,
intertwining operator algebras and the construction 
of genus-zero chiral conformal field theories. 
The main technical result used is the associativity 
for intertwining operators proved in \cite{H1} and \cite{H5}
for vertex operator algebras satisfying the three conditions 
stated in the theorem.

As is discussed in Section 2, assuming
the convergence and extension property and some other conditions,
the associativity for intertwining operators 
was proved in \cite{H1}.   In \cite{H3}, the commutativity for 
intertwining operators was proved using the associativity. 
In Lemma 4.1 in \cite{H6}, using the associativity and commutativity
proved in \cite{H1} and \cite{H2}, respectively,
it was shown  that one can analytically extend
products or iterates of two intertwining operators to 
$$M^{2}=\{(z_{1}, z_{2})\in \C^{2}\;|\; z_{1}, z_{2}\ne 0, z_{1}\ne
z_{2}\}.$$ 
Moreover, the construction of 
the map $(\Psi^{(1)}_{P(z_{1}, z_{2})})^{-1}$ in (14.65) in \cite{H1} 
can be reinterpreted as a proof of the fact that 
any four point correlation function can be obtained in this way.
For multipoint
correlation functions, the corresponding results can be proved using
the same method or using the results for four point correlation
functions above. This result shows that assuming the convergence 
and extension property and some other conditions, 
products or iterates of
intertwining operators can be analytically extended to multivalued
analytic functions defined on the moduli spaces of pointed
genus-zero Riemann surfaces.

The convergence and extension property was proved in \cite{H1}
under conditions weaker than the conditions we assume above.
The proof is based on the existence of systems of differential 
equations and the regularity of the singular points of 
these systems. Since we know only the existence, not the 
explicit forms, of these systems
of differential equations, many powerful tools available for 
Knizhnik-Zamolodchikov equations in the Wess-Zumino-Novikov-Witten
models are not available here anymore. 
However, the theory of intertwining operator 
algebras developed in \cite{H3}, \cite{H4} and \cite{H5}
allows us to reduce the use of these differential equations to a minimum.
For example, we do need the regularity of the singular points of 
these systems. But only the regularity of some ordinary 
differential equations induced from these systems are needed
and such regularity can be proved easily.

Now using the associativity for 
intertwining operators  repeatedly
to express the correlation functions obtained 
from products of 
three suitable intertwining operators as linear 
combinations of the correlation functions obtained 
from iterates of 
three intertwining operators in two ways, we obtain a formula 
for the matrix elements of the fusing isomorphisms. Then using 
some properties of the matrix elements of the fusing 
isomorphisms and their inverses proved in 
\cite{H7}, we obtain the first formula 
(\ref{formula1-2}).

The proof of the second formula (\ref{formula2}) 
mainly uses the results obtained in 
\cite{H6} on the convergence and analytic extensions 
of the $q_{\tau}$-traces of products of what we call ``geometrically-modified 
intertwining operators'', the genus-one associativity, and the modular
invariance of the space of these analytic extensions 
of the $q_{\tau}$-traces, where, as in Section 3,
$q_{\tau}=e^{2\pi \tau}$.
In \cite{Z}, in addition to the convergence and 
modular invariance of $q_{\tau}$-traces of vertex operators, 
Zhu also proved the convergence and modular invariance of the space of the 
$q_{\tau}$-traces of products of more than one 
vertex operators acting on modules. In \cite{DLM},
Dong, Li and Mason generalized Zhu's result to the case of twisted
modules. In \cite{M}, Miyamoto generalized Zhu's result to the case of
products of one intertwining operator and arbitrarily many vertex
operators.  However, these results of Zhu, Dong-Li-Mason and Miyamoto do
not give all genus-one correlation functions and thus are not 
enough for the construction
conformal field theories.  More specifically, the problem of
proving the convergence, the modular invariance and duality properties
for $q_{\tau}$-traces of products of {\it more than one} intertwining operator
was still open at that time. This problem is equivalent to the problem of 
constructing all the
chiral genus-zero correlation functions and establishing all the
desired properties.

The difficulty in the case of more than one intertwining operator is
that the method of Zhu, further developed by Dong-Li-Mason and
Miyamoto, cannot be adapted or generalized to this case, because there
is no commutator formula for general intertwining operators.  Here by
commutator formulas for general intertwining operators, we mean a
formula for 
$$\Y_{1}(w_{1}, x_{1})\Y_{2}(w_{2}, x_{2})-\Y_{3}(w_{2},
x_{2})\Y_{4}(w_{1}, x_{1})$$
where $\Y_{1}$, $\Y_{2}$, $\Y_{3}$ and $\Y_{4}$ are suitable
intertwining operators. There is no such formula, even in the case of
abelian intertwining operator algebras in the sense of Dong and Lepowsky
\cite{DL}.  Without such commutator formula, one cannot prove a
recurrence formula needed in the method of Zhu, Dong-Li-Mason and
Miyamoto.  Even the generalized commutator formula of Dong and Lepowsky
for abelian intertwining operator algebras cannot be used to obtained
such recurrence relations.  

In \cite{H6}, the author solved these previously open
problems.  These problems actually form one of the main obstructions to
a new range of applications of the theory of vertex operator algebras.
It is the solution to these problems together with the 
results on the duality properties of genus-zero correlation functions
discussed in the preceding section that allowed us to prove 
(\ref{formula2}). The method in \cite{H6} is completely different from
the one in \cite{Z}, \cite{DLM} and \cite{M} in the case of 
products of more than one operators: Instead of the nonexistent 
commutator formula, the author used the associativity and
commutativity for intertwining operators and the method 
of analytic extensions. The lack of a commutator formula was one of
the subtle difficulties which was overcome in \cite{H6}. 

To prove the second formula (\ref{formula2}), using the results 
proved in \cite{H6},
we introduce two operators $\alpha$ and $\beta$ on the space of 
linear maps from $\coprod_{a\in \A}W^{a}\otimes (W^{a})'$
to the space of genus-one two-point correlation functions
obtained from the analytic extensions of the 
$q_{\tau}$-traces of iterates of suitable ``geometrically-modified 
intertwining operators'' introduced in \cite{H6}:
$\alpha$ is induced from the 
translation of one of the points by $-1$ and $\beta$ is induced from
the 
translation of one of the points by $\tau$. The work in 
\cite{H6} is crucial in proving that these operators 
are well-defined and have the desired properties. The 
matrix $S$ also acts on the same space of linear maps
and it is easy to prove the relation
$$S\alpha S^{-1}=\beta.$$
We then use the genus-one associativity and other properties 
of the genus-one correlation functions in \cite{H6} to 
calculate $\alpha$ and $\beta$ explicitly in terms of 
the matrix elements of the fusing and braiding isomorphisms.
Substituting the explicit expressions of $\alpha$ and $\beta$
into this relation, we obtain (\ref{formula2}).

As in \cite{MS1}, the conclusions of the theorem follow immediately from
(\ref{formula1-2}) and (\ref{formula2}).

All the consequences of the Verlinde conjecture derived by physicists
now hold for vertex operator algebras satisfying the conditions in the 
theorem. In particular, we have the following Verlinde formula 
for fusion rules: For $a\in \mathcal{A}$, $S_{e}^{a}\ne 0$ 
and
\begin{equation}\label{v-form}
N_{a_{1}a_{2}}^{a_{3}}=
\sum_{a_{4}\in \A}\frac{S_{a_{1}}^{a_{4}}S_{a_{2}}^{a_{4}}S_{a_{4}}^{a'_{3}}}
{S_{e}^{a_{4}}}
\end{equation}
(cf. Section 3 in \cite{V}).
We also have the following well-known formula:
For $a_{1}, a_{2}\in \mathcal{A}$, 
\begin{equation}\label{s-form-3}
S_{a_{1}}^{a_{2}}
=\frac{S_{e}^{e}((B^{(-1)})^{2}(\Y_{a_{2}e; 1}^{a_{2}}
\otimes \Y_{a'_{1}a_{1}; 1}^{e};
\Y_{a_{2}e; 1}^{a_{2}}\otimes \Y_{a'_{1}a_{1}; 1}^{e}))}
{F(\Y_{a_{1}e; 1}^{a_{1}} \otimes \Y_{a'_{1}a_{1}; 1}^{e};
\Y_{ea_{1}; 1}^{a_{1}}\otimes \Y_{a_{1}a'_{1}; 1}^{e})
F(\Y_{a_{2}e; 1}^{a_{2}} \otimes \Y_{a'_{2}a_{2}; 1}^{e};
\Y_{ea_{2}; 1}^{a_{2}}\otimes \Y_{a_{2}a'_{2}; 1}^{e})}.
\end{equation}
The formula (\ref{s-form-3}) and certain properties of 
the matrix elements of the fusing and braiding isomorphisms
proved in \cite{H7} imply that the matrix $(S_{a_{1}}^{a_{2}})$ is 
symmetric.

\renewcommand{\theequation}{\thesection.\arabic{equation}}
\renewcommand{\thethm}{\thesection.\arabic{thm}}
\setcounter{equation}{0}
\setcounter{thm}{0}

\section{Rigidity, nondegeneracy and modular tensor categories}

A tensor category with tensor product bifunctor $\boxtimes$ 
and unit object $V$ 
is rigid if for every object 
$W$ in the category, there are right and left dual objects $W^{*}$
and $^{*}W$ together with morphisms $e_{W}: W^{*}\boxtimes W\to V$,
$i_{W}:V\to W\boxtimes W^{*}$, $e'_{W}: W\boxtimes ^{*}W\to V$
and $i'_{W}: V\to ^{*}W\boxtimes W$ 
such that the compositions of the morphisms in the sequence
$$\begin{CD}
W&@>>>
&V\boxtimes W
&@>i_{W}\boxtimes I_{W}>>
&(W\boxtimes W^{*})\boxtimes 
W&@>>>\\
&@>>>&W\boxtimes (W^{*}\boxtimes 
W)&@>I_{W}\boxtimes e_{W}>>&
W\boxtimes V&@>>>&W
\end{CD}$$
and three similar sequences are equal to the identity
$I_{W}$ on $W$. The rigidity is a standard notion in the theory 
of tensor categories. 
A rigid braided tensor category together with a twist (a natural 
isomorphism from the category to itself satisfying natural conditions)
is called a ribbon category. A semisimple ribbon category with 
finitely many inequivalent irreducible objects is a modular 
tensor category if the following nondegeneracy 
condition or  modularity is satisfied:
The $m\times m$ matrix formed by the traces
of the morphism $c_{W_{i}W_{j}}\circ c_{W_{j}W_{i}}$ 
in the ribbon category for 
irreducible modules $W_{1}, \dots, W_{m}$ is invertible. 
See \cite{T} and \cite{BK} for 
details of the notions in the theory of modular tensor 
categories. 

Using the results discussed in the proceeding section,
we obtain the following result:

\begin{thm}
Let $V$ be a simple vertex operator algebra satisfying the 
following conditions:

\begin{enumerate}

\item $V_{(n)}=0$ for $n<0$, $V_{(0)}=\mathbb{C}\mathbf{1}$, 
$W_{(0)}=0$ for any irreducible $V$-module which is 
not isomorphic to $V$.

\item Every $\mathbb{N}$-gradable weak $V$-module is completely 
reducible.

\item $V$ is $C_{2}$-cofinite, that is, $\dim V/C_{2}(V)<\infty$. 

\end{enumerate}
Then the braided tensor category structure on the 
category of $V$-modules constructed in 
\cite{HL1}--\cite{HL4}, \cite{H1} and \cite{H5}
is rigid, has a natural structure of 
ribbon category and satisfies the nondegeneracy condition. 
In particular, the category of $V$-modules 
has a natural structure of modular tensor category.
\end{thm}

Note that Condition 1 implies that
$V'$ is isomorphic to $V$ as a $V$-module. 
Thus Condition 1 in the theorem is slightly stronger
than Condition 1 in Theorem \ref{main}. 

We now discuss the proof of this theorem. 
We take both the left and right dual of a $V$-module $W$
to be the contragredient module $W'$ of $W$. Since our 
tensor category is semisimple, to prove the rigidity, we need 
only discuss irreducible modules. For $a\in \mathcal{A}$,
using the universal property
for the tensor product module $(W^{a})'\boxtimes W^{a}$,
we know that there exists a unique module map $\hat{e}_{a}:
(W^{a})'\boxtimes W^{a}\to V$ such that 
$$\overline{\hat{e}}_{a}(w'_{a}\boxtimes w_{a})
=\Y_{a'a; 1}^{e}(w'_{a}, 1)w_{a}$$
for $w_{a}\in W^{a}$ and $w'_{a}\in (W^{a})'$, where 
$\overline{\hat{e}}_{a}: \overline{(W^{a})'\boxtimes W^{a}}\to 
\overline{V}$ is the natural extension of $\hat{e}_{a}$. 
Similarly, we have a module map from
$W^{a}\boxtimes (W^{a})'$ to $V$. Since 
$W^{a}\boxtimes (W^{a})'$ is completely reducible and 
the fusion rule $N_{W^{a}(W^{a})'}^{V}$ is $1$, 
there is a $V$-submodule of $W^{a}\boxtimes (W^{a})'$
which is isomorphic to $V$ under the module map from 
$W^{a}\boxtimes (W^{a})'$ to $V$. 
Thus we obtain a module map
$i_{a}: V\to W^{a}\boxtimes (W^{a})'$ which maps
$V$ bijectively to this submodule of $W^{a}\boxtimes (W^{a})'$. 
Now
\begin{equation}\label{rigid-map}
\begin{CD}
W^{a}&@>>>
&V\boxtimes W^{a}
&@>i_{a}\boxtimes I_{W^{a}}>>
&(W^{a}\boxtimes (W^{a})')\boxtimes 
W^{a}&@>>>\\
&@>>>&W^{a}\boxtimes ((W^{a})'\boxtimes 
W^{a})&@>I_{W^{a}}\boxtimes \hat{e}_{a}>>&
W^{a}\boxtimes V&@>>>&W_{a}
\end{CD}
\end{equation}
is a module map from an irreducible module to itself. 
So it must be the identity map multiplied by a number.
One can calculate this number explicitly and it is equal to 
$$F(\Y_{ae; 1}^{a} \otimes \Y_{a'a; 1}^{e};
\Y_{ea; 1}^{a}\otimes \Y_{aa'; 1}^{e}).$$
From Theorem \ref{main}, this number is not $0$. 
Let 
$$e_{a}=\frac{1}{F(\Y_{ae; 1}^{a} \otimes \Y_{a'a; 1}^{e};
\Y_{ea; 1}^{a}\otimes \Y_{aa'; 1}^{e})}\hat{e}_{a}.$$
Then the map obtained from (\ref{rigid-map}) by 
replacing $\hat{e}_{a}$ by $e_{a}$ 
is the identity. 
Similarly, we can prove that all the other maps 
in the definition of rigidity are also equal to the 
identity. So the tensor category is rigid.

For any $a\in \A$, we define the twist on $W^{a}$ 
to $e^{2\pi i h_{a}}$. Then it is easy to verify that 
the rigid braided tensor category with this twist is 
a ribbon category.

To prove the nondegeneracy, we use the formula 
(\ref{s-form-3}). Now it is easy to calculate 
in the tensor category
the trace of $c_{W^{a_{2}}, W^{a_{1}}}\circ 
c_{W^{a_{1}}, W^{a_{2}}}$ for $a_{1}, a_{2}\in \mathcal{A}$, 
where $c_{W^{a_{1}}, W^{a_{2}}}: W^{a_{1}}\boxtimes 
W^{a_{2}}\to W^{a_{2}}\boxtimes 
W^{a_{1}}$ is the braiding 
isomorphism. The result is 
$$\frac{((B^{(-1)})^{2}(\Y_{a_{2}e; 1}^{a_{2}}
\otimes \Y_{a'_{1}a_{1}; 1}^{e};
\Y_{a_{2}e; 1}^{a_{2}}\otimes \Y_{a'_{1}a_{1}; 1}^{e}))}
{F(\Y_{a_{1}e; 1}^{a_{1}} \otimes \Y_{a'_{1}a_{1}; 1}^{e};
\Y_{ea_{1}; 1}^{a_{1}}\otimes \Y_{a_{1}a'_{1}; 1}^{e})
F(\Y_{a_{2}e; 1}^{a_{2}} \otimes \Y_{a'_{2}a_{2}; 1}^{e};
\Y_{ea_{2}; 1}^{a_{2}}\otimes \Y_{a_{2}a'_{2}; 1}^{e})}.
$$
By (\ref{s-form-3}), this is equal to 
$$\frac{S_{a_{1}}^{a_{2}}}{S_{e}^{e}}$$
which form an invertible matrix. So the semisimple
balanced rigid tensor category is nondegenerate.
Thus the tensor category is 
modular.  Details 
will be given in  \cite{H8}.

\noindent {\small \sc Department of Mathematics, Rutgers University,
110 Frelinghuysen Rd., Piscataway, NJ 08854-8019}

\noindent {\em E-mail address}: yzhuang@math.rutgers.edu

\end{document}